\RequirePackage{fix-cm}
\documentclass[smallextended]{svjour3}
\usepackage{mathptmx}
\usepackage{graphicx}
\usepackage[utf8]{inputenc}
\usepackage[numbers]{natbib} 
\usepackage{mathtools}
\usepackage{longtable}
\usepackage{listings}
\usepackage{booktabs}
\usepackage{algpseudocode}
\usepackage{amssymb}
\usepackage{algorithm}
\usepackage{multirow}
\usepackage{makecell}
\usepackage{booktabs}
\usepackage{biolinum}
\usepackage{subcaption}
\usepackage{caption}
\usepackage{hyperref}
\hypersetup{colorlinks, linkcolor=blue, citecolor=blue}
\usepackage[T1]{fontenc}
\usepackage{microtype}

\usepackage{url}

\usepackage{colortbl}
\usepackage{array}
\usepackage{graphicx}
\usepackage[table,xcdraw,svgnames,dvipsnames]{xcolor}
\usepackage{orcidlink}

\definecolor{bennevis}{HTML}{6baed6} 
\definecolor{almost_bennevis}{HTML}{c6dbef} 
\definecolor{benmacdui}{HTML}{fd8d3c} 
\definecolor{cairngorm}{HTML}{fdd0a2} 
\definecolor{wider_bennevis}{HTML}{74c476} 
\definecolor{top_25_munros}{HTML}{a1d99b} 
\definecolor{top_50_munros}{HTML}{756bb1} 
\definecolor{top_135_munros}{HTML}{bcbddc} 
\definecolor{mountainous}{HTML}{636363} 
\definecolor{lowland}{HTML}{bdbdbd} 

\usepackage{etoolbox} 
\newcommand\fig[2][]{\hyperref[fig:#2]{Figure~\ref*{fig:#2}\ifstrempty{#1}{}{#1}}}

\algnewcommand\Not{\textbf{not} }
\newcommand{\optimaNum}{957,174\ }

\title{Where's Ben Nevis?}
\subtitle{A 2D optimisation benchmark with \optimaNum local optima based on Great Britain terrain data}

\author{
Yuhang Wei \orcidlink{0009-0007-8822-690X} \and 
Michael Clerx \orcidlink{0000-0003-4062-3061} \and 
Gary R. Mirams \orcidlink{0000-0002-4569-4312}
}

\institute{Yuhang Wei \and Michael Clerx \and Gary R. Mirams \at School of Mathematical Sciences, University of Nottingham, Nottingham, NG7 2RD, UK \\ Corresponding author:  Gary R. Mirams \\ \email{gary.mirams@nottingham.ac.uk}}

\begin{document}

\maketitle


\abstract{
	We present a novel optimisation benchmark based on the real landscape of Great Britain (GB). 
	The elevation data from the UK Ordnance Survey Terrain 50 dataset is slightly modified and linearly interpolated to produce a target function that simulates the GB terrain, packaged in a new Python module \texttt{nevis}. 
	We introduce a discrete approach to classifying local optima and their corresponding basins of attraction, identifying \optimaNum local optima of the target function. 
	We then develop a benchmarking framework for optimisation methods based on this target function, where we propose a Generalised Expected Running Time performance measure to enable meaningful comparisons even when algorithms do not achieve successful runs (find Ben Nevis). 
	Hyperparameter tuning is managed using the \texttt{optuna} framework, and plots and animations are produced to visualise algorithm performance. 
	Using the proposed framework, we benchmark six optimisation algorithms implemented by common Python modules. 
 Amongst those tested, the Differential Evolution algorithm implemented by \texttt{scipy} is the most effective for navigating the complex GB landscape and finding the summit of Ben Nevis.
}

\section{Introduction}


%
%

Optimisation benchmarks are crucial tools for evaluating and comparing the performance of various optimisation algorithms \citep{bartzbeielstein2020benchmarking}.
Numerous 2D optimisation benchmarks exist, many of which are detailed in sources such as the Wikipedia page on test functions for optimisation\footnote{ \url{https://en.wikipedia.org/wiki/Test_functions_for_optimisation\#Test_functions_for_single-objective_optimisation}}. 
These benchmarks are often called `artificial landscapes', carefully designed mathematical functions that simulate complex optimisation problems.

However, with the availability of Ordnance Survey open terrain elevation data, we saw an opportunity to develop a benchmark problem based on the \textit{real landscape} of the entirety of Great Britain (GB). 
This real-world landscape presents a unique and challenging optimisation testbed, significantly different from artificial benchmarks typically used in optimisation studies.


Our motivations for developing this benchmark problem are multi-faceted:

\begin{itemize}
	\item \textbf{Educational Tool:} By using the real landscape of Great Britain, we can create an educational tool that is intuitive and engaging.
	This tool can help students and practitioners visualise and understand how different optimisation algorithms navigate complex landscapes, offering clear insights into their workings and performance.
	\item \textbf{Challenging Local Optima:} The landscape of Great Britain is characterised by a large number of local optima. 
	In our study, we identified \optimaNum local optima based on a discrete mathematical definition applied to the Ordnance Survey 2D grid points, which approximately simulates the real-world concept of mountain peaks, as explained in Section~\ref{sec:local-optima}. 
	This significant number of local optima raises a fundamental question: Can any optimisation method reliably find the global optimum in such a complex landscape? 
	Investigating this will provide valuable insights into the strengths and limitations of various optimisation techniques.

	\item \textbf{Suitability of Optimisers:} By analysing where different optimisers converge in the GB landscape, we aim to learn more about their suitability for different types of problems.
	Understanding the final positions and paths of various optimisers can guide us in selecting appropriate algorithms for specific optimisation challenges, thereby enhancing their practical applicability.
\end{itemize}

In summary, the real landscape of Great Britain offers a rich and challenging benchmark for optimisation algorithms, providing a novel and practical alternative to traditional artificial landscapes. 
Through this study, we hope to contribute to the understanding and development of more effective optimisation methods.

\subsection{Great Britain elevation data}

\begin{figure}
	\centering
	\includegraphics[width=0.85\linewidth]{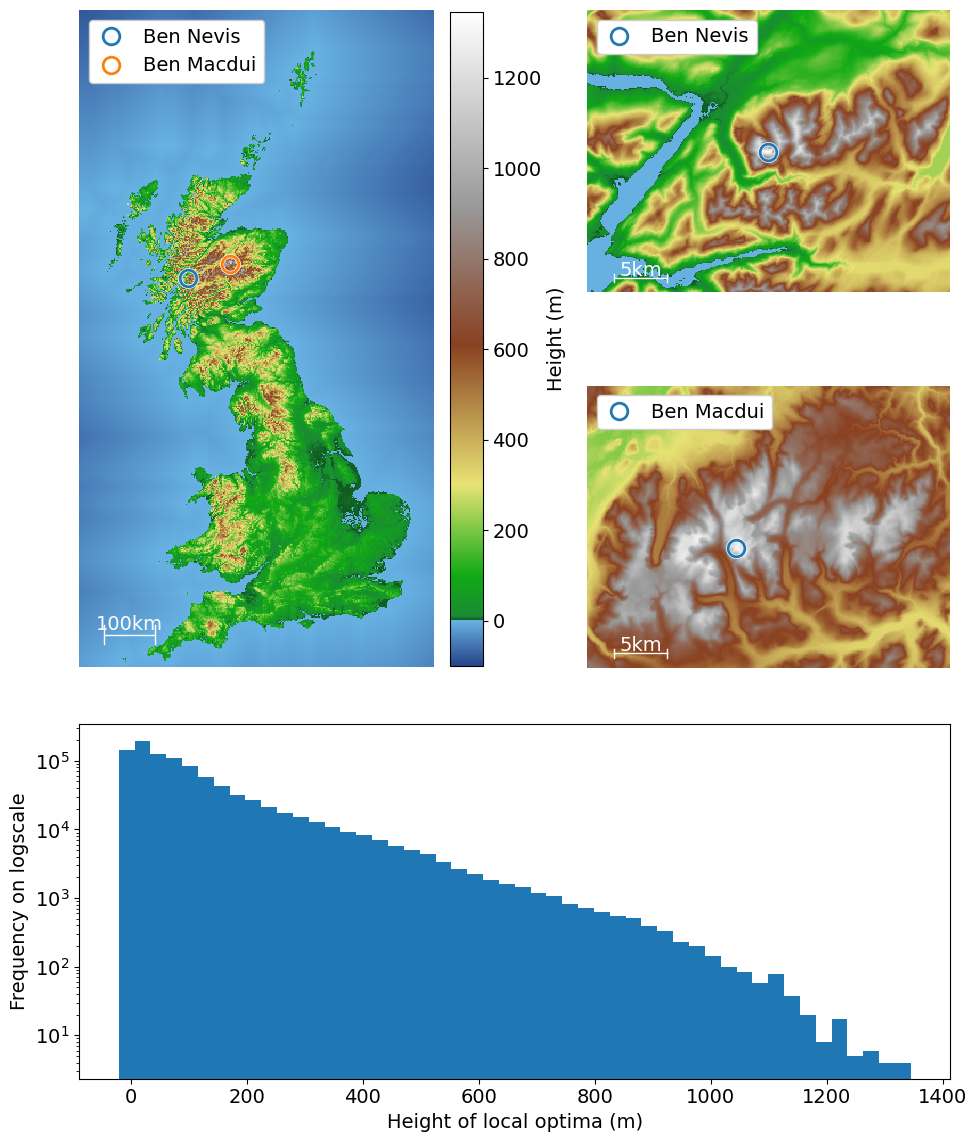}
	\caption{Great Britain terrain data adapted from the Ordnance Survey `OS Terrain 50' dataset.
    A 700km by 1300km rectangle is shown, including the land mass of Great Britain and its islands, but omitting nearby landmasses Ireland and France, and with algorithmically generated seabed elevation levels (see main text).
    Insets show the area around the two highest peaks (Ben Nevis 1st and Ben Macdui 2nd) on the same scale, note how Ben Macdui is part of a far larger area of high peaks and plateaus.
   The lower panel displays the logarithmic frequency distribution of the heights of all \optimaNum local optima}
	\label{fig:map}
\end{figure}

\fig{map} shows a map of the terrain data, and Table~\ref{tab:height_bands} defines some height bounds (objective function values) that are particularly of interest for Great Britain.
A little geopolitics and hill walking terminology is in order at this point.
Great Britain (GB) consists of a mainland (occupied by the bulk of England, Scotland and Wales) and outlying islands such as the Shetlands, Orkneys, Hebrides, Isle of Wight, Isles of Scilly \& Anglesey.
GB is not the same as the United Kingdom (UK), as the UK also includes Northern Ireland.
GB also excludes self-governing islands such as the Isle of Man and the Channel Islands.

Exactly what constitutes a distinct mountain/hill in the fractal nature of the landscape has been defined in the past with somewhat subjective criteria.
The highest mountains in GB are all in Scotland, and an early attempt at a definition and ranking of all Scottish mountains over 3,000\,ft (914.4m) was made by Sir Hugh Munro in 1891 \citep{munro1891}. 
The mountains on his tables became known as `The Munros', with adjoining sub-peaks also over 3,000\,ft known as `Munro Tops'.
There are now 282 Munros on the list that is maintained by the Scottish Mountaineering Club (revisions over time have occurred due to general attempts at rationalisation and more recently modern Global Positioning System survey methods), alongside a list of `compleators'\footnote{which stands at 7831 people as of September 2024 (\url{https://smc.org.uk/hills/compleators})}.

Of all British peaks (and hence Munros) Ben Nevis is the highest at 1345\,m, it is situated on the West coast of the Scottish Highlands and is surrounded on three sides by terrain going back to almost sea level, but its immediate massif to the East includes the 7th--9th highest Munros (see zoom region in Figure~\ref{fig:map}).
The second highest Munro is Ben Macdui, approximately 85\,km east in the Cairngorm mountain range, itself in close proximity to 3rd--6th highest Munros such that an exhausting day's walk could include all of these.
Compared to Ben Nevis, the Cairngorms generally form a much more considerable high plateau of arctic weather, flora and fauna (see second zoom region, on the same scale, in Figure~\ref{fig:map}).
This will become relevant when looking at where various optimisers tend to terminate.
Six mountains in England and fifteen in Wales also meet the 3,000\,ft threshold, and these are known as `Furths'.
For a complete description of various lists of GB hills and their definitions, we refer the reader to the excellent Database of British and Irish Hills\footnote{\url{http://www.hills-database.co.uk/database_notes.html\#defs}}.

\begin{table}[tbh]
	\centering
	\caption{Bands of height displayed in other figures and scores used in hyperparameter optimisation
    }
	\label{tab:height_bands}
    \begin{tabular}{@{}l p{3.6cm} c c c l@{}}
        \toprule
        Interval (m) & Optima Description & Score & Basin Size & Basin Proportion & Colour \\
        \midrule
        \([-100, 0)\) & Below sea-level (fens) & --- & \(1.23 \times 10^{5}\) & \(3.37 \times 10^{-4}\) & \cellcolor{lowland} \\
        \([0, 600)\) & Lowland areas & --- & \(3.58 \times 10^{8}\) & \(9.83 \times 10^{-1}\) & \cellcolor{lowland} \\
        \([600, 1000)\) & Mountainous areas & --- & \(5.77 \times 10^{6}\) & \(1.58 \times 10^{-2}\) & \cellcolor{mountainous} \\
        \([1000, 1100)\) & Mountains within top \(\sim\)135 Munros \& 5 Welsh Furths & --- & \(2.61 \times 10^{5}\) & \(7.18 \times 10^{-4}\) & \cellcolor{top_135_munros} \\
        \([1100, 1150)\) & Mountains within top \(\sim\)50 Munros & --- & \(6.65 \times 10^{4}\) & \(1.83 \times 10^{-4}\) & \cellcolor{top_50_munros} \\
        \([1150, 1215)\) & Mountains within top \(\sim\)25 Munros & --- & \(4.14 \times 10^{4}\) & \(1.14 \times 10^{-4}\) & \cellcolor{top_25_munros} \\
        \([1215, 1235)\) & Wider Ben Nevis Massif height (within top 9 Munros) & 1 & \(4.99 \times 10^{3}\) & \(1.37 \times 10^{-5}\) & \cellcolor{wider_bennevis} \\
        \([1235, 1297)\) & Cairngorm Plateau height (within top 6 Munros) & 2 & \(1.18 \times 10^{4}\) & \(3.24 \times 10^{-5}\) & \cellcolor{cairngorm} \\
        \([1297, 1310)\) & Ben Macdui (2nd highest Munro) & 3 & \(1.32 \times 10^{3}\) & \(3.63 \times 10^{-6}\) & \cellcolor{benmacdui} \\
        \([1310, 1340)\) & On Ben Nevis (not quite at summit) & 7 & \(2.18 \times 10^{3}\) & \(5.99 \times 10^{-6}\) & \cellcolor{almost_bennevis} \\
        \([1340, 1346)\) & Ben Nevis (highest Munro) & 10 & \(8.78 \times 10^{2}\) & \(2.41 \times 10^{-6}\) & \cellcolor{bennevis} \\
        \toprule
    \end{tabular}
\end{table}
\section{Methods}

\subsection{Data and modifications}

To create a dataset for optimisation, we downloaded and slightly modified the `OS Terrain 50' dataset, made available by Great Britain's Ordnance Survey (OS)\footnote{\url{https://www.ordnancesurvey.co.uk/products/os-terrain-50}} under an Open Government Licence (OGL).
This contains both contour data, which was not used in this study, and a rectangular grid of points with given elevation (in metres, with a vertical resolution of 0.1\,m.) spaced 50\,m apart. 
The data are clustered in 10\,km by 10\,km `tiles', and tiles are provided to cover the entire landmass of GB and its islands, resulting in a 700\,km by 1300\,km rectangle as shown in \fig{map}.
No data are provided for the elevation of the seabed, or indeed of the landmasses of Ireland and France which would otherwise have been visible within this area.

To account for the missing data, we implemented an artificial sloping seabed, such that points further away from coastlines are increasingly far below sea level.
This approach guides the optimisers towards the land, which is the area of interest --- a completely flat objective of 0\,m in the sea would lead to many optimisers becoming `stuck' and would not teach us much about their behaviour on the GB landscape.
Because the full data set, encoded as a matrix of single-precision floating point numbers, is 1.6\,GB in size, this modification required a somewhat complicated divide-and-conquer approach, which is detailed in Appendix~\ref{app:sea}.

\subsubsection{Interpolation}
The data are presented to the optimisers through a callable function $h(x, y)$ that takes a grid coordinate (in meters) as input, and returns a height (again in meters).
Grid coordinates are relative to the lower-left (most southwestern) grid point, and are known as `Eastings' ($x$) and `Northings' ($y$).
Because we are interested in testing optimisation on a continuous domain, $h$ implements a linear interpolation on the grid.
While not the most realistic, a linear interpolation is fast, does not require the storage of interpolation coefficients, and does not introduce artificial maxima (as is prone to happen with e.g. cubic splines).

\subsubsection{Packaging and distribution}
To facilitate use of the modified OS Terrain 50 dataset for optimisation benchmarking, we created a fully documented Python module called \texttt{nevis}.
Upon first use, users are prompted to download the compressed data set directly from the Ordnance Survey.
Once downloaded, the sea-slope generating code is run, and results are stored to disk.
In addition to providing the interpolating function $h$, \texttt{nevis} provides plotting functionality (as showcased in \fig{map}) and conversion from grid coordinates to online mapping systems such as Google Maps.
The \texttt{nevis} module supports Python 3.6 and newer, is published on GitHub\footnote{\url{https://github.com/CardiacModelling/BenNevis}} and available via the Python Package Index (PyPI)\footnote{\url{https://pypi.org/project/nevis/}}.

\subsection{Local optima and basins of attractions}
\label{sec:local-optima}

In the context of the OS Terrain 50 data, we characterise the \textit{neighbours} of a given grid point \(x = (i, j)\) as the set of 8 grid points surrounding our point (including diagonally-adjacent), that is the set of \(x' = (i', j')\) that satisfy the condition \(\max \{ |i - i'|, |j - j'| \} = 1\).
These eight neighbours are given the following `canonical' anti-clockwise ordering: bottom, bottom-right, right, top-right, top, top-left, left, bottom-left. 


We define the \textit{gradient} from grid point \((i, j)\) to \((i', j')\) as
\begin{equation}
	g((i, j), (i', j')) = \frac{h(i', j') - h(i, j)}{\sqrt{(i - i')^2 + (j - j')^2}}.
\end{equation}

We attempt to assign to any given grid point \((i, j)\) its  \textit{neighbour of steepest-ascent (NSA)} according to the following steps:
\begin{enumerate}
    \item If there exists a neighbour \((i', j')\) whose height is strictly greater than that of \((i, j)\), then the NSA of \((i, j)\) is the neighbour \((i', j')\) with the maximum gradient \(g((i, j), (i', j'))\). 
    If two or more neighbours share the same maximum gradient, we choose the one which appears the first in the canonical ordering.
    \item We now iterate and try to assign an NSA to any grid point \((i, j)\) which has not been assigned one yet (which can only have happened if all the neighbours are of equal height or lower). 
    In the \(k\)-th iteration, for each such \((i,j)\), we try to find a neighbour \((i', j')\) with the same height as \((i, j)\) that has been assigned an NSA.
    Then if \((i', j')\) has an \(l\)-th order NSA \((i'', j'')\) higher than \((i', j')\) (and thus higher than \((i, j)\)), where \(l \le k\) is the smallest such integer, then the NSA of \((i, j)\) is the neighbour \((i', j')\) with maximal \(g((i, j), (i'', j''))\). 
    The \(l\)-th order NSA of any grid point is the point obtained by recursively applying the NSA operation \(l\) times on it.
    Again, any tie is resolved by resorting to the canonical ordering. 
    If in the last iteration, there is no point whose NSA gets updated, then the algorithm has converged and no more iterations will be carried out.
    \item Finally, we would like to ensure that there is at most one local optimum in a connected block of points with the same height. 
    This is achieved by launching breadth-first searches (BFS) on all grid points currently without an NSA.
    If the current step of the BFS is at grid point \((i, j)\), then we expand to any neighbour \((i', j')\) with the same height and also without an NSA, and  \((i', j')\) obtains its NSA as \((i, j)\).
    \item Now any grid point \((i, j)\) that still does not have an NSA is defined as a \textit{local optimum}.
\end{enumerate}

Using this approach we identified \optimaNum local optima. 

There are 282 Munros and 227 Tops, giving a total of 509 Scottish summits over 3,000 ft in the current Munro Tables.
But the OS data contain many more local optima which are not distinct enough to be generally defined as separate mountains, but which do form local optima on the surface we wish to optimise over.
Of our \optimaNum optima, 1,242 in GB have a height over 3,000\,ft (914.4\,m), and out of these 1,192 are in Scotland, giving a sense of how many more optima we are working with relative to what a mountaineer might typically consider a separate summit.

It is not hard to show that the steps above effectively implement a partial ordering on the grid points if we consider each grid point as `less than' its NSA (if its NSA exists), and there is an underlying directed acyclic graph  \(G = (V, E)\) induced by this partial ordering, where the vertex set \(V\) is the set of all grid points, and
a directed edge from grid point \(x\) to \(x'\) is in the edge set \(E\) if and only if \(x'\) is the NSA of \(x\).
A local optimum is exactly a sink vertex of the graph \(G\), i.e. a vertex with out-degree zero.

With this understanding, we establish a discretised version of the gradient ascent algorithm, which we refer to as algorithm \(\cal{A}\).
Starting from any grid point, we iteratively progress to the NSA of the current point, until we reach a local optimum, where the algorithm terminates.

Let \(x_0\) be any given local optimum.
We define the \textit{basin of attraction (BoA)} of \(x_0\) as the collection of initial points that converge to \(x_0\) via algorithm \(\cal{A}\).
In other words, a point \(x\) is contained in the BoA of \(x_0\) if and only if there exists a path from \(x\) to \(x_0\) on graph \(G\).
The \textit{area} of the BoA of \(x_0\) is the number of grid points it contains.

To effectively ascertain the BoA and its area for each local optimum, we could use  breadth-first searches by starting the search from each local optimum on the reverse graph of \(G\).
Alternatively, we could carry out the following procedure to label each grid point with the BoA it belongs to: 
\begin{enumerate}
    \item Label each local optimum with the BoA of itself;
    \item Start from any grid point which has not been labelled yet, and iteratively progress to the NSA of the current point whilst recording all the points which we visit, until we reach a point whose BoA has been labelled. Then label all the recorded points with the BoA;
    \item Repeat step 2, until all grid points are labelled.
\end{enumerate}

Table \ref{tab:height_bands} lists \(11\) height bands and their corresponding descriptions. We refer to the `below sea-level' interval as \(H_0\), `lowland areas' as \(H_1\), and so on.
We say that the BoA of a local optimum \(x_0\) belongs to a height band \(H_j = [h_j, k_j)\) if \(x_0\) lies in \(H_j\). 
The \textit{basin size} of \(H_j\) is the area sum of all BoAs that belong to \(H_j\), and the \textit{basin proportion} of \(H_j\) is the proportion of its basin size in the total number of grid points. 
The data are shown in the last two columns of Table \ref{tab:height_bands}. 
The basin proportion of \(H_j\) indicates the probability for a single run of algorithm \(\cal{A}\) with a (uniformly) random starting point to obtain a local optimum in \(H_j\).
For example, the probability of a single run of algorithm \(\cal{A}\) with a random starting point to result in the Ben Nevis height band is estimated to be  \(2.41 \times 10^{-6}\).
We highlight that a local optimum identified through this discrete process aligns with a local optimum on the linear interpolant surface, which serves as the objective function of our benchmarking framework.
However, the concept of `basin of attractions' is dependent on the local optimiser used; in particular, the basin of attraction can be defined as the set of initial points that lead to the respective local optimum when applied to that specific optimiser.
Hence, our definition of the `basin of attractions' which relies on algorithm \(\cal{A}\) may not coincide with the `basin of attractions' defined using any given continuous optimiser that is more commonly used.

\subsection{Algorithms}

\subsubsection{Basic concepts}

\begin{definition}
	An ``algorithm'' (or an optimiser) is a fixed procedure to find the optimum of the objective function, with or without a hyperparameter space.
\end{definition}

\begin{definition}
	A ``hyperparameter'' of an algorithm is a value that the algorithm makes use of and that needs to be assigned before the algorithm is run. The hyperparameter space of an algorithm is the range of values that its hyperparameters can take.
\end{definition}

\begin{definition}
	An ``algorithm instance'' is an algorithm with all of its hyperparameters assigned to a particular value in the hyperparameter space.
\end{definition}

We qualitatively distinguish between global and local optimisers:
\begin{itemize}
	\item A \textit{global optimiser} is designed to find the global optimum of the objective function, which is the best possible solution across the entire search space. 
	Global optimisers aim to avoid getting trapped in local optima by employing strategies that allow for extensive exploration of the search space. 
	\item  A \textit{local optimiser} focuses on finding a local optimum, which is the best solution within a limited region of the search space. 
	Local optimisers are often more efficient in terms of computation but are prone to getting trapped in local optima if the search space is complex and rugged.
\end{itemize}

We remark that sometimes it is not clear-cut how to classify an algorithm as either global or local, because this distinction can be blurred, especially when an algorithm employs techniques from both categories to balance exploration and exploitation effectively or when its effectiveness depends heavily on specific problem characteristics and hyperparameter settings.

One common method to extend the capability of a local optimiser to approximate a global search is the \textit{multistart} approach, where the local optimiser is run multiple times from random initial points distributed throughout the search space. 
Each run of the local optimiser aims to find a local optimum, and by aggregating the results of multiple runs, we increase the likelihood of identifying the global optimum.

On the other hand, many global optimisers employ a two-phase structure to enhance their performance \citep{floudas_encyclopedia_2009}:
\begin{itemize}
	\item The \textit{global phase} focuses on exploring the search space broadly to identify promising regions. 
	Techniques such as random sampling, evolutionary strategies, or other global search methods are used to cover a wide area and avoid premature convergence to local optima.
	\item After identifying promising regions in the global phase, the algorithm switches to a \textit{local phase} to refine the optimisation results. Local search methods, such as gradient-based algorithms or other optimisers like Nelder--Mead, introduced below, are used to fine-tune the solutions within these regions.
\end{itemize}


\subsubsection{List of algorithms used}

The algorithms that we use are listed along with their hyperparameter space in Table~\ref{tab:algo-list}.
Here are brief descriptions of the algorithms and their implementations:
\begin{itemize}
     \item \textbf{Nelder--Mead} \citep{nelder_simplex_1965} is a derivative-free optimisation algorithm designed for nonlinear, multidimensional unconstrained problems. 
     It operates using a simplex, a geometric figure formed by \( n+1 \) vertices in an \( n \)-dimensional space, which iteratively adjusts to converge towards the minimum of the objective function. 
     The algorithm adapts the simplex through operations such as reflection, expansion, contraction, and shrinkage to efficiently explore the solution space and refine the search around the minimum. 
     The Nelder--Mead implementation we use is provided by the NLOPT library \citep{nlopt2021neldermead}, which offers a collection of algorithms suitable for a wide range of optimisation problems. 
     The library is available in multiple programming languages, including C, C++, Python, and MATLAB, and we use its Python interface.
     Note that there is no hyperparameter for Nelder--Mead and the multistart version of it is used as a global optimiser.
    \item \textbf{Dual Annealing} \citep{xiang_generalized_1997} is a stochastic global optimisation algorithm combining the concept of simulated annealing with local search strategies.  
    Dual Annealing operates by iteratively exploring the solution space at decreasing ``temperatures'', which allows it to accept worse solutions at higher temperatures to escape local minima. 
    The local search component refines the solutions found during the annealing process, improving convergence efficiency. 
    In the Python ecosystem, Dual Annealing is implemented in the SciPy library \citep{scipy2024dual}, a widely-used library for scientific computing.
    We choose Nelder--Mead as the local search strategy for Dual Annealing.
    \item \textbf{Multi-Level Single-Linkage (MLSL)} \citep{rinnooy1987stochastic1, rinnooy1987stochastic2} is a two-phase global optimisation algorithm that combines a global search strategy with a local optimisation method. MLSL begins by sampling the solution space and applying clustering to group the samples. 
    It then performs local optimisation within each cluster, iteratively refining the clusters to focus on more promising areas. This hierarchical approach attempts to balance exploration and exploitation.
    The MLSL implementation we use is by NLOPT \citep{nlopt2021mlsl}, with Nelder--Mead being the local optimisation method.
    \item \textbf{CMA-ES (Covariance Matrix Adaptation Evolution Strategy)} \citep{hansen_cma_2016} is an evolutionary algorithm for nonlinear, non-convex optimisation problems. 
    It operates by iteratively sampling candidate solutions, evaluating their fitness, and updating the covariance matrix to adapt the search distribution towards promising regions of the solution space. 
    CMA-ES is designed to handle complex optimisation landscapes by effectively balancing exploration and exploitation. 
    We use the CMA-ES implementation by the PINTS package \citep{clerx2018probabilistic}.
    \item \textbf{Particle Swarm Optimisation (PSO)} \citep{kennedy_particle_1995} is a population-based stochastic optimisation technique inspired by the social behaviour of birds flocking or fish schooling. 
    PSO operates by initialising a swarm of particles in the search space, where each particle represents a candidate solution. 
    These particles move through the search space, adjusting their positions based on their own experience and the experience of neighbouring particles. 
    The algorithm iteratively updates the velocity and position of each particle to explore the solution space and converge towards the optimal solution. 
    We use the PSO algorithm as implemented in the PINTS package \citep{clerx2018probabilistic}.
    \item \textbf{Differential Evolution} \citep{storn_differential_1997} is a population-based global optimisation algorithm particularly well-suited for continuous spaces. 
    Differential Evolution operates by iteratively improving a candidate solution with regard to a given measure of quality or fitness. 
    It initialises a population of potential solutions and uses operations such as mutation, crossover, and selection to evolve the population over generations. 
    The algorithm relies on the differences between randomly sampled pairs of solutions to explore the search space and guide the population towards the global optimum. 
    We use the Differential Evolution algorithm as implemented in the SciPy library \citep{scipy2024de}.
\end{itemize}
A few remarks on the hyperparameters listed in Table~\ref{tab:algo-list}:
\begin{itemize}
    \item All hyperparameters are documented in detail by the corresponding implementation of each algorithm.
    \item For Differential Evolution, there is a hyperparameter \texttt{mutation} which is derived as follows: if \texttt{dithering} is set to \texttt{True}, then \texttt{mutation} is the interval \([\texttt{mutation\_low}, \texttt{mutation\_high}]\), otherwise \texttt{mutation} is the value \texttt{mutation\_low}.
\end{itemize}

\begin{table}[tbh]
	\centering
	\caption{Algorithms with their hyperparameter spaces}
	\label{tab:algo-list}
	\begin{tabular}{p{0.18\textwidth}p{0.25\textwidth}p{0.1\textwidth}p{0.25\textwidth}p{0.1\textwidth}}

		\toprule
		Name & Hyperparameters               & Default value & Range                                  & Type
		\\ \midrule

		\multirowcell{4}{ Dual Annealing                                                                               \\ (\texttt{scipy}) }
		     & \texttt{initial\_temp}        &      \(5230\)                 & \([0.2, 5\times10^4]\)                   & Log    \\
		     & \texttt{restart\_temp\_ratio} &        \(2 \times 10^{-5}\)               & \([1\times10^{-6}, 0.9]\)                & Log    \\
		     & \texttt{visit}               &       \(2.62\)                & \([1.5, 2.9]\)                          & Linear \\
		     & \texttt{accept}                &        \(-5\)               & \([-5, -1.1 \times 10^{-4}]\)            & Linear \\

		\midrule

		\makecell{MLSL                                                                                                 \\ (\texttt{nlopt)}  }
		     & \texttt{population}           &          \(4\)             & \([1, 3000]\)                            & Int    \\

		\midrule

		\multirowcell{2}{ CMA-ES \\ (\texttt{PINTS}) }
		     & \texttt{sigma0}               &           \(1.2 \times 10^{5}\)            & \( [3.5 \times 10^4, 3.5 \times 10^5] \) & Linear \\
		     & \texttt{population\_size}       &          \(6\)             & \([4, 100]\)                               & Int    \\
            \midrule

            \multirowcell{3}{PSO  \\ (\texttt{PINTS}) }
                & \texttt{sigma0}               &           \(1.2 \times 10^{5}\)            & \( [3.5 \times 10^4, 3.5 \times 10^5] \) & Linear \\
		     & \texttt{r}       &          \(0.5\)             & \([0,1]\)                               & Linear    \\
		     & \texttt{population\_size}             &          \(6\)             & \([1, 1000]\)                              & Int    \\

            \midrule

            \multirowcell{6}{Differential \\ Evolution \\ (\texttt{scipy})}
            & \texttt{popsize}             &          \(15\)             & \([10, 50]\)                              & Int    \\
            & \texttt{recombination}    &   \(0.7\)         &\([0, 1]\)         & Linear \\
            & \texttt{polish} & \texttt{False} &  & Boolean \\
            & \texttt{dithering} & \texttt{True} &  & Boolean \\
            & \texttt{mutation\_low} & 0.5 & \([0,2]\) & Linear \\
            & \texttt{mutation\_high} & 1 & \([\texttt{mutation\_low}, 2]\) & Linear \\
		\bottomrule
	\end{tabular}
\end{table}


As we will see, the hyperparameters of a given algorithm can significantly impact their performance on our test function and therefore hyperparameter tuning methods will be employed, which in effect introduces an additional layer of optimisation. 
However, while we strive for improved performance, ensuring hyperparameters are `optimally' configured may not be feasible or necessary. 
Instead, our focus is on reaching a balance where the configurations are adequately tuned to yield enhanced results within the constraints of our computational resources and time limitations.

\subsection{Benchmarking}

\subsubsection{Objective}
The objective of our benchmarking process is to evaluate specific algorithm implementations that are widely used, with regard to their performance on the GB terrain data (height function).
Modifications, such as adjustments to hyperparameters, are made only when the interface permits.

\subsubsection{Termination criteria}
Our primary aim is to maintain as much consistency as possible, so we apply the same termination criteria across various algorithms.
However, it is important to note that the same settings may be treated differently by various implementations, which may influence the benchmarking results.

An algorithm is terminated when one of the following \textit{global criteria} is applicable and met:
\begin{enumerate}
	\item When it reaches a target function value \(f_\text{target}\) close to the height of Ben Nevis; or
	\item When the number of function evaluations reaches a budget \(T_\text{max}\).
\end{enumerate}

Further, the following \textit{local criteria} are applied to a local optimiser or a local phase of a global optimiser:
\begin{enumerate}
	\item When an optimisation step changes the function value by less than the absolute tolerance \(f_\text{tol}\); or
	\item When an optimisation step changes the value of each parameter by less than the absolute tolerance \(x_\text{tol}\).
\end{enumerate}
Moreover, any algorithm may terminate on its own by design apart from being subject to these criteria. 

Usually, a global algorithm is subject to the global criteria, while its local phases are subject to the local criteria. 
However, as mentioned earlier, it is not always clear-cut whether an algorithm should be classified as global or local. 
To address this ambiguity, we introduce the concept of an \textit{early-terminating algorithm}, defined as an algorithm that may terminate for reasons other than the global criteria.

To ensure consistency across our evaluations, we apply a multi-start strategy for all algorithms. 
This means that the multi-start version of each algorithm is subject to the global criteria, ensuring that these criteria are the only reasons for termination in any run, thereby standardising the termination process.
For non-early-terminating algorithms, this approach does not change their behaviour. 
The details of this strategy are discussed further in Section \ref{sec:multi-start}.

We observe that certain algorithms may not fully comply with these criteria.
For example, in some algorithms, the number of function evaluations during a run may slightly exceed the budget \(T_\text{max}\).
This may occur because, for example, the algorithm does not support termination during its local phase or always evaluates its whole internal population for particle-based optimisers.
To address this issue, we employ a strategy in which all function calls made during a run are recorded.
Subsequently, we consider either the initial \(T_\text{max}\) function calls or those up to the first occurrence (inclusive) where the function value reaches or surpasses the threshold \(f_\text{target}\), depending on whichever condition is met first.
From this subset of function calls, we identify the highest function value attained, which we refer to as the \textit{returned height} of the run.

We also note that the implementation of \(f_\text{tol}\) and \(x_\text{tol}\) may exhibit some variation across different implementations.
For instance, in the \texttt{nlopt} implementation of the Nelder--Mead algorithm, the algorithm terminates when either the absolute tolerance \(f_\text{tol}\) or the absolute tolerance \(x_\text{tol}\) is met.
In contrast, the \texttt{scipy} implementation of the same algorithm requires both tolerances to be satisfied for termination. As mentioned before, we do not interfere with these differences but only keep the values set consistently.
Table~\ref{tab:termination_criteria} shows the values we set to each termination criterion in our experiment.

\begin{table}[ht]
	\centering
	\caption{Prescribed values in optimisation termination criteria. For reference, the height of Ben Nevis is 1345m}
	\label{tab:termination_criteria}
	\begin{tabular}{ll}
		\toprule
		Variable            & Value    \\
		\midrule
		\(f_{\text{target}}\) & 1,340\,m \\
		\(T_{\text{max}}\)    & 50,000   \\
		\(f_{\text{tol}}\)    & 0.2\,m   \\
		\(x_{\text{tol}}\)    & 10\,m    \\
		\bottomrule
	\end{tabular}
\end{table}

\subsubsection{Initial guesses}

In order to guarantee reproducibility, for each run of an algorithm instance, the initial guess is pseudorandomly generated based on the run index. 
Upon initializing NumPy's random number generator with a random seed, which is equal to the run index, we generate two pseudorandom numbers within $[0, 1)$ from a uniform distribution \citep{NumPyDocs}. 
These numbers are scaled linearly to the maximum dimensions of the rectangular domain of the target function, with width $7\times 10^5$ and height $1.3 \times 10^6$. 

\subsubsection{Performance measures}
\label{sec:performance}
Consider a stochastic algorithm instance that has run \(N\) times with different random seeds and initial guesses using the termination criteria above.
For each run, the number of function evaluations and returned height are recorded.
A \textit{performance measure} for this instance over the \(N\) runs is a scalar value obtained by aggregating the run results.
Generally, two objectives can be used as the criteria to measure the performance of algorithms: consumed budget and solution quality, which give rise to the two following questions:
\begin{itemize}
	\item Horizontal-cut view: How fast can algorithms reach a given solution quality, i.e.\ a predefined target function value \(f_\text{target}\)?
	\item Vertical-cut view: What solution quality can the algorithms achieve within a given number of function evaluations?
\end{itemize}
The names `horizontal-cut' and `vertical-cut' view come from drawing a horizontal or vertical line in a convergence graph (see Section~\ref{sec:plots}).
In our discussion of performance measures, our primary focus is on the horizontal-cut view.
This perspective offers better interpretability, allowing us to assert that an algorithm instance is a certain number of times faster than another.
The vertical-cut view is used in performance plots.

We classify a run of an algorithm instance as \textit{successful} if its returned height is larger than \(f_\text{target}\); otherwise, it is \textit{unsuccessful}.

Denote the index set of all \(N\) runs as \(I = \{  1, \dots, N \}\).
The index set of successful runs is denoted as \(I_s\), and let \(N_s = |I_s|\) be the number of successful runs.
Let \(T_i\) denote the number of function evaluations of the \(i\)-th run. 
Note that by our termination criteria and by applying the multi-start strategy to early-terminating algorithms, if the \(i\)-th run is unsuccessful, then \(T_i = T_\text{max}\).
We can then define the \textit{success rate} of an algorithm instance over \(N\) runs as
\begin{equation}
	p_s = \frac{N_s}{N}.
\end{equation}

A commonly used performance measure is the \textit{expected running time} (ERT), defined as the expected number of function evaluations needed for the algorithm to reach \(f_\text{target}\) once:
\begin{equation}
	\text{ERT} = \frac{\sum_{i \in I} T_i} {N_s}.\label{eq:ert}
\end{equation}
So a lower ERT indicates better performance. 
When \(N_s = 0\), i.e.\ there is no successful run, we define $\mathrm{ERT} = \infty$.
When unsuccessful runs all use the full \(T_\text{max}\) evaluations that are allowed,  ERT has an alternative expression \citep{auger_performance_2005},
\begin{equation}
	\text{ERT} = T_s + T_{\text{max}}\frac{1 - p_s }{p_s},
\end{equation}
where we define the \textit{mean running time of successful runs}
\begin{equation}
	T_s = \frac{1}{N_s} \sum_{i \in I_s} T_i.
\end{equation}

Therefore, an algorithm instance with a low ERT would potentially have: i) a high \(p_s\) (the algorithm should be able to find Ben Nevis or other high Munros --- see Table~\ref{tab:height_bands} --- as frequently as possible); and ii) a low \(T_s\) (when the algorithm successfully does so, it should do so within a low budget) \citep{bossek_multi-objective_2020}.

However, using ERT directly in our context can be problematic if an algorithm almost never successfully finds Ben Nevis before its termination.
Some algorithms, such as CMA-ES, tend to be `trapped' in large basins of attraction such as the Cairngorms.
Even though such algorithms could have a high chance of discovering the second highest Munro, Ben Macdui, we frequently get an infinity value of the ERT of their instances, which makes comparing instances and hyperparameter tuning impractical.
As a result, we introduce the following performance measure, as a generalisation of the ERT:
\begin{equation}
	\text{GERT} = \frac{\sum_{i \in I} T_i} {\sum_{i \in I} S_i},
\end{equation}
where \(S_i\) is the \textit{score} of the \(i\)-th run, assigned according to the interval in which the returned height falls in Table~\ref{tab:height_bands}.
GERT (Generalised ERT) can be interpreted as the expected number of function evaluations to win one such score point.
We can see that ERT represents a specific case of GERT, where the interval \([1340, 1345)\) (measured in meters) is assigned a score of \(1\), while all other intervals receive scores of \(0\).

We also mention a performance measure, \textit{averaged returned height}, the average of the returned heights of the \(N\) runs, which is based on the vertical-cut view. 
Some further performance measures are listed in Appendix~\ref{app:further-pm}.

\subsubsection{Multi-start strategy}
\label{sec:multi-start}

Let \(\cal{E}\) be an early-terminating algorithm with \(N\) runs.
For simplicity, assume that the last run is successful.
To calculate the performance measures of \(\cal{E}\), we need to implement \(\cal{E}\) as its multi-start version, which we call algorithm \(\cal{E}'\).
Now the \(N\) early-terminating runs of algorithm \(\cal{E}\) may be merged into \(N'\) multi-start runs of algorithm \(\cal{E}'\) for \(N' \le N\).
(The last run of \(\cal{E}\) is successful, then so is the last merged run of \(\cal{E}'\); in this way we ensure that the last run of \(\cal{E}'\) is not early-terminating.)
We are interested in the differences in the performance measures calculated based on the early-terminating runs and the multi-start runs. 

ERT is almost invariant under this multi-start transformation;
strictly speaking, the only difference here is that in a multi-start run, the last call to algorithm \(\cal{E}\) may be truncated when \(T_{\text{max}}\) is reached.
This in general should only have a small influence on both the denominator and numerator of Equation~\ref{eq:ert}.

However, other performance measures, including GERT, success rate and averaged returned height, are not invariant under this transformation.
As an example, suppose algorithm \(\cal{E}\) easily discovers the height range above 1215m but rarely discovers Ben Nevis.
Further suppose that \(\cal{E}\) terminates on its own after exactly $500$ function evaluations in each run and it is run $100$ times, where it discovers the height range $[1235, 1297)$ in $20$ runs, discovers Ben Nevis exactly at the last function evaluation of the $100$-th run, and does not gain scores otherwise. 
With multi-start, these $100$ runs are $100$ local phases of a single global run, which gains $10$ scores for the discovery of Ben Nevis in the end; 
however, without multi-start, each run can gain scores on its own, giving  $20 \times 2 + 10 = 50$ scores gained in total. 
This reasoning indicates that, without multi-start, early-terminating algorithms gain an unfair advantage in terms of GERT over global algorithms.
The difference between ERT and GERT is caused by the fact that the termination criteria only rely on the discovery of Ben Nevis and not on gaining scores otherwise. 
Thus, the multi-start approach is necessary to ensure the consistency of the performance measures when applied to algorithms with different termination patterns.




\subsubsection{Hyperparameter tuning}
\label{sec:ht}

\begin{algorithm}[hbt!]
	\caption{Hyperparameter tuning based on GERT with pruning}\label{alg:hpt}
	\begin{algorithmic}
		\For{\(k \in \left\{ 0, 1, \dots, M-1\right\}\) }
		\State \(I_k \gets\) a new algorithm instance from the sampling algorithm
		\State \(N \gets 0\)
		\While{\(\sum_{i=0}^{N-1} T_i < T_\text{instance}\)}
		\State \(R_N \gets \) run result of \(I_k\)
		\State \(T_N \gets\) the number of function evaluations of \(R_N\)
		\State \(\text{GERT}_N \gets \) the GERT of \(I_k\) over \(R_0, \dots, R_{N}\)
		\State report \(\text{GERT}_N\) to the pruning algorithm
		\If{the pruning algorithm decides \(I_k\) is not promising }
		\State prune \(I_k\), move on to the next instance
		\EndIf
		\State \(N \gets N + 1\)
		\EndWhile
		\State \(I_k.\text{GERT} \gets \text{GERT}_{N-1}\)
		\If{\(k=0\) or \(I_k.\text{GERT} < I_\text{best}.\text{GERT}\)}
		\State \(I_\text{best} \gets I_k\)
		\EndIf
		\EndFor
	\end{algorithmic}
\end{algorithm}

Due to the large scale of the problem (indicated by the large number of local optima and the small Ben Nevis basin of attraction), most algorithms do not perform well with their default hyperparameter settings (if applicable).
We thus need hyperparameter tuning so that the algorithms can `explore' the landscape sufficiently to produce reasonably good results.
\texttt{optuna} \citep{akiba2019optuna}, an automatic hyperparameter tuning framework commonly used in machine learning, is utilised for this purpose.

The hyperparameter tuning process, shown in Algorithm~\ref{alg:hpt}, can be outlined as follows:
A sequential generation of \(M\) algorithm instances is executed, employing a hyperparameter space sampling algorithm.
This sampling algorithm takes into account the hyperparameters and GERT of all prior instances with an attempt to ascertain their correlations and explore the hyperparameter space efficiently.
Each instance is allocated a fixed total number of function evaluations, denoted as \(T_\text{instance}\), distributable across individual runs.
(In practice, \(T_\text{instance}\) can be exceeded by the last run; precisely, \(T_\text{instance} \le \sum_{i=1}^{N} T_i < T_\text{instance} + T_\text{max}\).)
Notably, during each run, the instance's current GERT is computed and reported to the pruning algorithm, which may prematurely terminate an instance if it is deemed unpromising.
Finally, the instance with the lowest GERT is selected as the optimal representation of the algorithm under consideration. 
We employ the default sampling and pruning algorithms of \texttt{optuna}.
Prescribed values used in the process can be seen in Table~\ref{tab:hpt-values}.

\begin{table}[ht]
	\centering
	\caption{Prescribed values in hyperparameter tuning}
	\label{tab:hpt-values}
	\begin{tabular}{ll}
		\toprule
		Variable              & Value     \\
		\midrule
		\(M\)                   & 100       \\
		\(T_{\text{instance}}\) & 1,000,000 \\
		\bottomrule
	\end{tabular}
\end{table}

To ensure both the meaningfulness of comparisons between different algorithm instances and the reproducibility of our experiments, we consistently assign the same random seed and initial guess to the \(i\)-th run of each instance.
However, we remark that most global optimisation algorithms have an internal sampling method, which we do not wish to interfere with.
This means that the initial guess may not play a significant role in the whole optimisation process, probably unless the initial guess is close enough to the global optimum.
Also, since different sampling methods may treat the same random seed differently, it may not be meaningful to compare the results produced by different algorithms of a single run even with the same random seed.

For each instance, we fix the total number of function evaluations \(T_\text{instance}\) instead of the number of runs \(N\).
This can make sure that early-terminating algorithms can be measured in a larger number of runs.


\subsubsection{Performance plots} \label{sec:plots}


Consider an algorithm instance with \(N\) runs such that the \(i\)-th run has \(T_{i}\) function evaluations for \(i = 1, \dotsc, N\) and \(1 \le T_i \le T_\text{max}\). 
The performance plots try to answer the question: how does the algorithm instance perform at a certain number of function evaluations based on these runs?

Consider first for simplicity the case when \(T_\text{max} = T_1 = T_2 = \dotsc = T_N\). For the \(i\)-th run, we can find the largest height \(F_{i,j}\) among the first \(j\) function evaluations for each \(1 \le j \le T_i = T_\text{max}\). 
Then \(F = \{  F_{i,j} \}_{(i, j)}\) can be thought of as a matrix with \(N\) rows and \(T_\text{max}\) columns.
The \(i\)-th row of matrix \(F\) shows how the `solution quality' of the \(i\)-th run increases with the increase of the number of function evaluations, corresponding to the `horizontal-cut view', whereas the \(j\)-th column of matrix \(F\) reveals the `solution quality' of the algorithm instance given a budget of \(j\) function evaluations based on the \(N\) runs, corresponding to the `vertical-cut view'. 
We can then calculate the mean and five-number summary (minimum, first quartile, median, third quartile, and maximum) of the \(j\)-th column and show them in a line graph, and this is called the \textit{aggregated convergence graph}.
On the other hand, we can count the number of entries of the \(j\)-th column that lie within each of the height band Table~\ref{tab:height_bands}, and show them in a stacked graph, which we call the \textit{height-band graph}.

Now due to our termination criteria, when a run is successful, it would in fact have function evaluations less than \(T_\text{max}\). 
In this case, the `matrix' \(F\) would effectively have some empty entries at the end of each row.
A simple remedy is to `pad' each row by filling these empty entries with the last available number of each row.
This is sensible as these padded numbers are still the best-so-far heights of the run in the first \(j\) function evaluations.
Also, since the padding procedure only takes place for successful runs, the padded value would only be within the Ben Nevis height band (Table~\ref{tab:height_bands}).
This padding technique is applied in plotting the aggregated convergence graphs and height-band graphs. 


We also produce a variant of the aggregated convergence graph, where \(F_{i, j}\) is replaced by the closest distance to Ben Nevis among the first \(j\) function evaluations for the \(i\)-th run, with everything else the same as above.

In addition, for any algorithm instance with \(N\) runs, we could calculate its ERT with different \(f_\text{target}\) values, for any \(f_\text{target} \le 1340\). 
For such a given \(f_\text{target}\), we can re-classify the \(i\)-th run as successful or not by checking if there is function evaluation with height at least \(f_\text{target}\) its first \(T_\text{max}\) function evaluations. 
If yes, then the run is successful and \(T_i\) is the number of function evaluations prior to this point (inclusive). 
Otherwise, the run is unsuccessful and \(T_i = T_\text{max}\).
In this way, we can calculate the ERT value based on the given \(f_\text{target}\) using Equation~\ref{eq:ert}. 
This leads to a line graph that shows the correlation between ERT and \(f_\text{target}\) for each algorithm instance.




\subsubsection{Animations}

To better study the behaviour of the optimisers, we produced animations for each run of an algorithm instance. 
These animations illustrate the evaluated points in each iteration of the run against the backdrop of the entire Great Britain map and a regional map around Ben Nevis.
The visual representation allows us to observe the dynamic search patterns and convergence behaviours of the different algorithms.

\subsection{Computational Environment}

The experiments were conducted on an Ubuntu 22.04.3 LTS system running on the Windows Subsystem for Linux (WSL2) environment, with a kernel version \texttt{5.15.153.1-microsoft-standard-WSL2}. 
The hardware configuration includes an 11th Gen Intel Core i7-1165G7 processor, which has 4 cores and 8 threads, with a base clock speed of 2.80 GHz. 
The system is equipped with 7.6 GiB of RAM, with approximately 5.9 GiB available. 

Python version 3.10.12 was used for executing the code with the following Python package versions: \texttt{nevis} (0.1.0), \texttt{scipy} (1.14.1), \texttt{nlopt} (2.8.0), \texttt{pints} (0.5.0), \texttt{optuna} (4.0.0) and \texttt{numpy} (2.1.1).
In particular, all plots and visualisations were generated using the Python package \texttt{matplotlib} (3.9.2).
The code for finding local optima and their basins of attraction requires C++17, compiled with \texttt{g++} (11.4.0), and the \texttt{pybind11} module (2.13.6) for integration with Python.
A complete list of packages and their versions is provided in the accompanying repository.

\section{Results}

Note that the areas of the basins of attraction (Table~\ref{tab:height_bands}), and their proportion of the total area, provide a useful guide as to how many function evaluations might be needed. 
The basin of attraction for Ben Nevis covers a proportion of only $2.41\times10^{-6}$, suggesting that a naive strategy of picking an initial guess at random, and doing gradient descent algorithms would only work for one in $1/(2.41\times10^{-6})=414,938$ attempts. Each of these attempts would need a number of function evaluations to complete the gradient descent, so the expected number of function evaluations for this strategy would be $\mathcal{O}(10^6)$ to $\mathcal{O}(10^7)$.

\begin{figure}
    \centering
    \includegraphics[width=\linewidth]{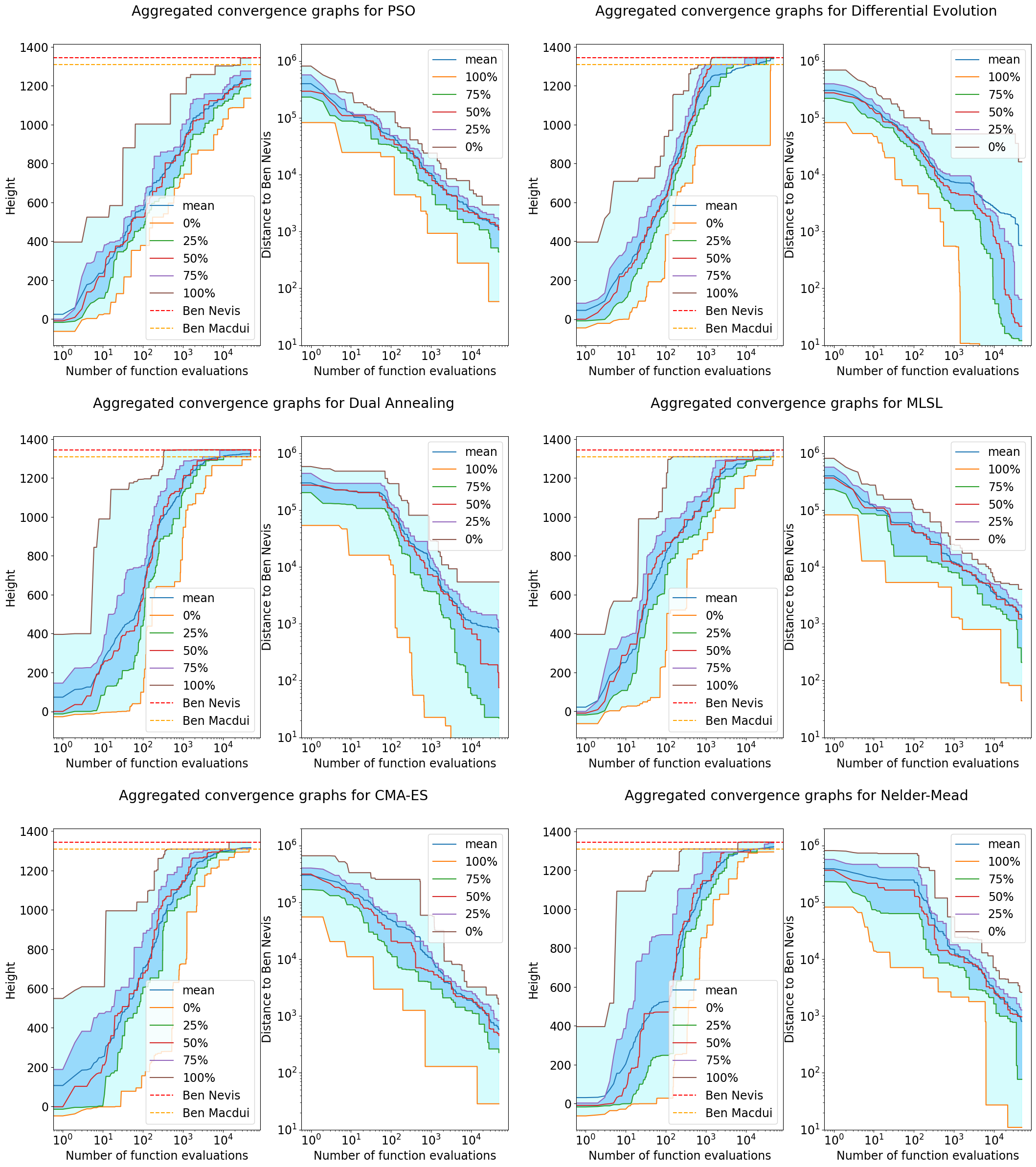}
\caption{Aggregated convergence graphs for the selected instance of each algorithm. 
	In each aggregated convergence graph on the left, the \(x\)-axis represents the number of function evaluations, and the lines show the mean and five-number summary of the best-so-far heights across all runs at each function evaluation.
	For runs that terminate before the maximum number of function evaluations, the best-so-far height is `padded'. 
	Each aggregated convergence graph on the right shows the closest-so-far distance to Ben Nevis instead of the best-so-far height}
    \label{fig:agg-1}
\end{figure}

\begin{figure}
    \centering
    \includegraphics[width=\textwidth]{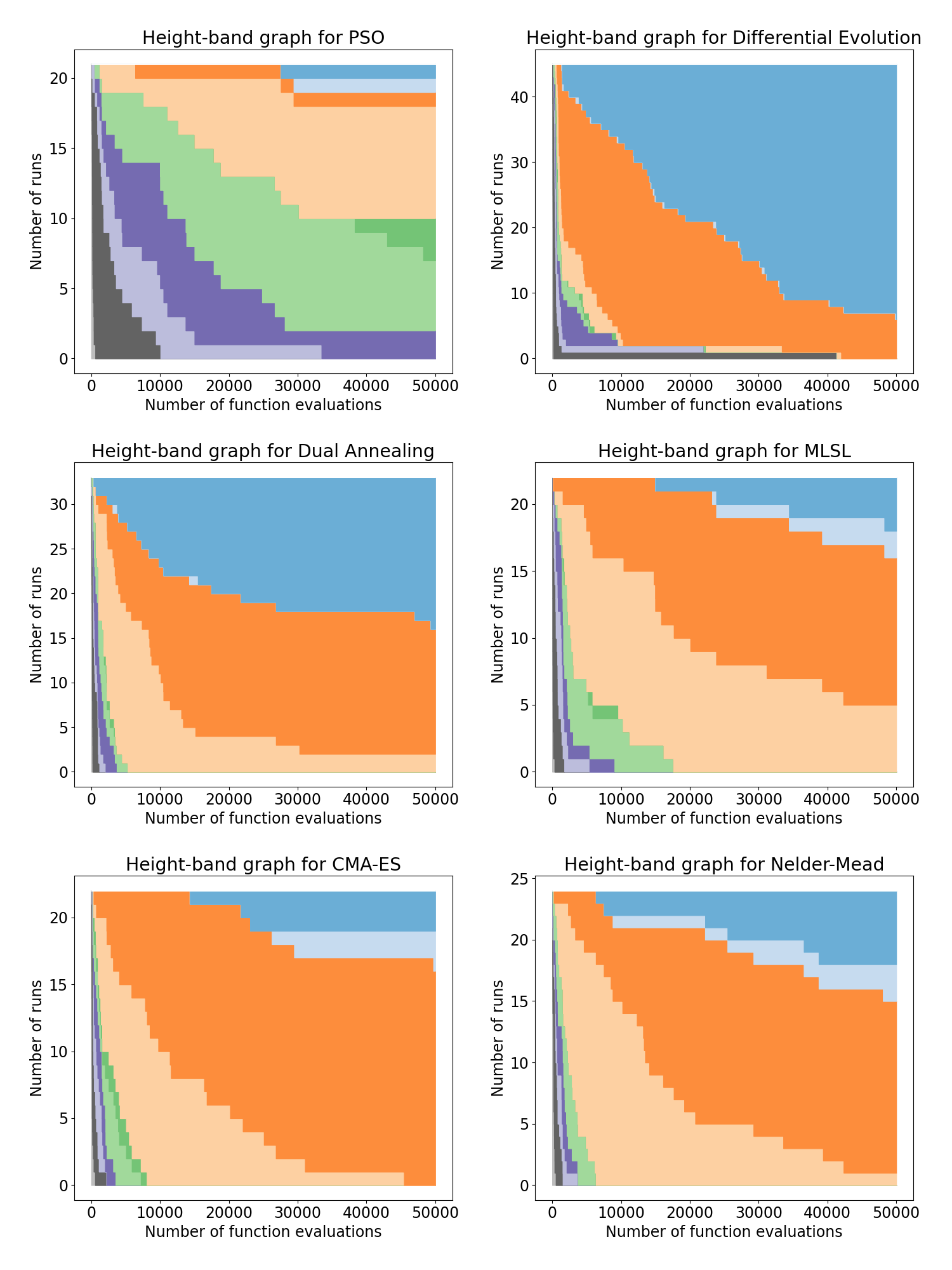}
	\caption{Height-band graphs for the selected instance of each algorithm. Each height-band graph displays the distribution of the best-so-far heights within predefined height bands (Table~\ref{tab:height_bands}) across all runs at each number of function evaluations. 
	For runs that terminate before the maximum number of function evaluations, the best-so-far height is `padded'. 
	The dark blue (Ben Nevis) region of each graph is the result of the padding process}
	\label{fig:hb-1}
\end{figure}

\begin{figure}
	\centering
	\includegraphics[width=0.6\linewidth]{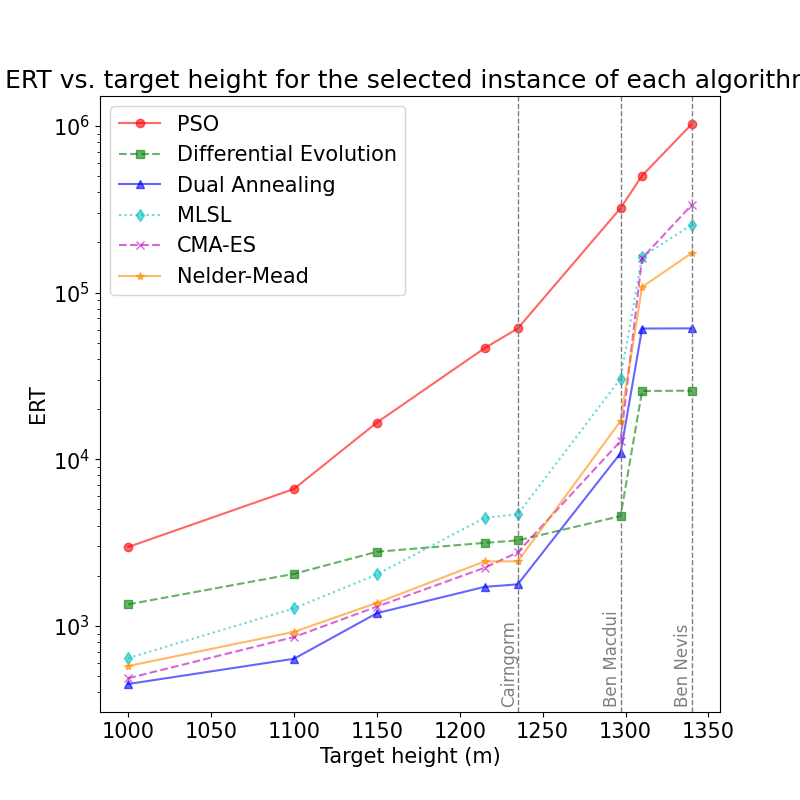}
	\caption{ 
		ERT with different target heights in Table~\ref{tab:height_bands} based on the runs of the selected instance of each algorithm. 
		The heights of Cairngorm Plateau, Ben Macdui and Ben Nevis are labelled as vertical lines
	}
	\label{fig:ert-plot}
\end{figure}

We evaluated the performance of six optimisation algorithms, Nelder--Mead, Dual Annealing, MLSL, CMA-ES, PSO, and Differential Evolution, using our proposed benchmarking framework.
The hyperparameters of the selected instance of each algorithm after tuning are shown in Table~\ref{tab:best-instance}, and the performance measures of each selected instance Table~\ref{tab:performance-measures}. 
The aggregated convergence graphs for each selected instance are plotted in Figures~\ref{fig:agg-1}, and the height-band graphs Figures~\ref{fig:hb-1}. 
The line graph for ERT with different \(f_\text{target}\) for all instances is shown in Figure~\ref{fig:ert-plot}.
Among the six algorithm instances, Nelder--Mead, PSO, CMA-ES and Differential Evolution tend to terminate early and their performance measures are based on their corresponding `multi-start' algorithm, as explained in Section~\ref{sec:multi-start}.

We first observe that the four performance measures give nearly identical rankings of the six instances, except that CMA-ES has better averaged returned height and GERT than MLSL but performs worse in terms of the other two measures.
We also remark that algorithms with a higher success rate have more run numbers since we have a fixed budget \(T_\text{instance}\) for all runs of an instance and a successful run terminates before reaching \(T_\text{max}\).

Differential Evolution outperformed the other algorithms in terms of all listed measures; notably, it demonstrated a success rate of 87\%, with the lowest GERT \(2.47 \times 10^3\). 
We observe that in each iteration, Differential Evolution tends to first sample points across the whole map, then concentrate sampling around Scotland or other promising areas, polish the results with a local optimiser, and terminate early, usually only with hundreds of function evaluations. 
This strategy ensures a full exploration of the landscape before the algorithm is committed to any local optima.

The second best algorithm is Dual Annealing, with a success rate slightly above 50\% and a GERT approximately twice that of Differential Evolution. 
Its annealing process helps the algorithm not get trapped in local optima and gives a preferable strategy for global exploration.
The relatively high performance of Dual Annealing may have also benefited from its high-dimensional hyperparameter space.

Interestingly, Nelder--Mead (multi-start version) performs better than MLSL with a Nelder--Mead local phase. 
This may be because the \texttt{nlopt} implementation of MLSL fixed a few hyperparameters\footnote{\url{https://github.com/stevengj/nlopt/blob/7cdebfe5f777b12d3c5b0788c38fe595444d69c6/src/algs/mlsl/mlsl.c\#L315-L325}}, making it impossible to fine-tune them for better performance on our landscape.
In addition, the global sampling strategy of MLSL, which tries to avoid repeated searches of the same local optima, may be prone to neglecting the Ben Nevis `basin of attraction', which is relatively small and close to low attitude regions such as the sea.
In other words, the complication of the GB landscape indicates that inappropriate \textit{a priori} assumptions that guide the global sampling process may be detrimental to the performance. 

The GERT and average returned height of CMA-ES indicate that its performance is on par with that of MLSL. 
The height-band graph for CMA-ES shows that the algorithm is highly effective at discovering Ben Macdui, the second highest Munro.
Furthermore, as shown in Figure~\ref{fig:ert-plot}, the ERT with Ben Macdui as the target height for CMA-ES ranks the third best, behind Differential Evolution but close to Dual Annealing.
This is likely because the Cairngorm Plateau region, with its high average elevation, is relatively easy to locate and favours the evolutionary process of CMA-ES.

Finally, PSO has the lowest success rate of 5\%, where only one successful run is conducted.
The plots reveal that PSO is not even very effective in discovering the Cairngorm Plateau. 
Since PSO updates particles' positions based on both their own best-known positions and the global best position, in a rugged landscape, particles might frequently change directions, leading to inefficient search patterns and slow convergence.

\begin{table}[tbh]
	\centering
	\caption{Best instances selected for each algorithm}
	\label{tab:best-instance}
	\begin{tabular}{p{0.2\textwidth}p{0.3\textwidth}p{0.25\textwidth}}

		\toprule
		Name & Hyperparameters               & Value
		\\ \midrule
		\multirowcell{4}{ Dual Annealing                             }
		& \texttt{initial\_temp} & \(2.69 \times 10^{4}\) \\
		& \texttt{restart\_temp\_ratio} & \(1.49 \times 10^{-3}\) \\
		& \texttt{visit} & \(2.47\) \\
		& \texttt{accept} & \(-3.42\) \\
		\midrule

		\makecell{MLSL                                                 }
		& \texttt{population} & 15 \\
		\midrule

		\multirowcell{2}{ CMA-ES }
		& \texttt{sigma0} & \(3.09 \times 10^{5}\) \\
		& \texttt{population\_size} & 10 \\

        \midrule
        \multirowcell{3}{PSO}
		& \texttt{sigma0} & \(1.94 \times 10^{5}\) \\
		& \texttt{r} & \(2.68 \times 10^{-1}\) \\
		& \texttt{population\_size} & 124 \\

 \midrule
 \multirowcell{6}{Differential \\ Evolution}
 & \texttt{popsize} & 11 \\
 & \texttt{recombination} & \(6.77 \times 10^{-1}\) \\
 & \texttt{polish} & True \\
 & \texttt{dithering} & True \\
 & \texttt{mutation\_low} & \(7.50 \times 10^{-1}\) \\
 & \texttt{mutation\_high} & \(9.18 \times 10^{-1}\) \\
		\bottomrule
	\end{tabular}

\end{table}

\begin{table}[tbh]
	\caption{Performance measures of the best instance of each algorithm}
	\label{tab:performance-measures}
	\centering
	\begin{tabular}{cccccc}
		\toprule
		Name                        & \makecell{Run                                                          \\number} & \makecell{Success\\rate} & \makecell{Average\\returned\\height (m)} & ERT & GERT \\
		\midrule
		Differential Evolution	& 45	& 87\%	& 1340	& \(2.58 \times 10^{4}\) 	& \(2.47 \times 10^{3}\) \\
				\midrule
		Dual Annealing	& 33	& 52\%	& 1326	& \(6.10 \times 10^{4}\) 	& \(4.80 \times 10^{3}\) \\
		\midrule
        Nelder--Mead	& 24	& 25\%	& 1320	& \(1.73 \times 10^{5}\) 	& \(8.30 \times 10^{3}\) \\

		\midrule
        CMA-ES	& 22	& 14\%	& 1316	& \(3.36 \times 10^{5}\) 	& \(1.02 \times 10^{4}\) \\
        \midrule
		MLSL	& 22	& 18\%	& 1313	& \(2.55 \times 10^{5}\) 	& \(1.05 \times 10^{4}\) \\
		\midrule
        PSO	& 21	& 5\%	& 1236	& \(1.03 \times 10^{6}\) 	& \(2.63 \times 10^{4}\) \\
		\bottomrule
	\end{tabular}
\end{table}

\section{Discussion}


We have proposed an optimisation benchmarking framework based on GB landscape data.
We developed the Python module \texttt{nevis} that facilitates the use of the modified OS Terrain 50 dataset for optimisation benchmarking. 
In the context of 2D grid points, we introduced an approach to defining and finding local optima and their corresponding basins of attraction in a way that resembles the classification of mountains, discovering \optimaNum local optima from the OS Terrain 50 data.

We proposed a new performance measure, GERT, as a generalisation of the existing ERT measure to facilitate the comparison of algorithm instances in cases where there are no successful runs. 
To address the potential influence of early termination on performance measures, we suggested a `multi-start' approach to ensure the consistency of these measures across algorithms with different termination patterns.
For hyperparameter tuning, we utilised the \texttt{optuna} framework, which supports efficient exploration in the hyperparameter space by pruning unpromising instances. 
We also produced visualisations of algorithm performance, including aggregated convergence graphs and height-band graphs, and demonstrated the optimisation process using animations.


Finally, we conducted a benchmarking process on six optimisation algorithms using implementations from mainstream Python modules, aiming to interpret the results. 
Our findings highlight Differential Evolution as the top-performing algorithm across all measures, achieving a success rate of 87\% and a GERT of \(2.47 \times 10^3\). 
The results suggest that algorithms prioritising global exploration—whether through effective sampling of promising regions, as seen with Differential Evolution, or the adaptive annealing process of Dual Annealing—tend to perform better in our landscape. 
Moreover, the comparison between Nelder–Mead and MLSL illustrates how inflexible hyperparameter settings or inappropriate global sampling processes can hinder the algorithm performance.
Additionally, CMA-ES appears more adept at identifying Munros on the Cairngorm Plateau, likely because the higher average elevation of the region is favoured by the evolution process of the algorithm. 
Overall, these findings indicate that algorithms emphasising global exploration and offers higher flexibility of its hyperparameter space tend to excel in our landscape.

Optimisation algorithms often work differently in higher dimensions, and the intuition we gain by understanding their behaviour for 2D problems does not always translate. 
It would be possible to create higher dimensional problems based on the GB terrain, for instance a 4D benchmark for parameters $i,j,k,l$ with a `height' defined by $h(i,j,k,l) = \sqrt{h(i,j) \times h(k,l)}$.

Future research can extend our benchmarking framework in several ways, such as incorporating additional real-world landscapes from different geographical regions or exploring smaller, localised areas within Great Britain to analyse specific optimisation challenges in greater detail. 
Importantly, our framework is designed for easy integration of new algorithms, allowing researchers to add emerging techniques for straightforward comparisons in terms of a range of performance measures and plots. 
This modularity not only enriches the analysis of algorithm performance in navigating complex terrains but also encourages ongoing innovation in optimisation strategies.


In conclusion, our benchmarking framework based on the real landscape of Great Britain has provided valuable insights into the performance of various optimisation algorithms. 
The identification of \optimaNum local optima highlights the complexity of real-world landscapes and the challenges they present for optimisers.
Differential Evolution (with `multi-start' strategy) emerged as the leading algorithm, emphasising the importance of global exploration in successfully navigating such terrains. 
Expanding this framework to incorporate additional landscapes, algorithms, and higher-dimensional problems could further enhance our understanding and application of optimisation techniques.



\bibliography{source.bib}

\section*{Statements and Declarations}

\subsection*{Funding}
This work was supported by the Wellcome Trust (grant no. 212203/Z/18/Z) via a Senior Research Fellowship to GRM.

\subsection*{Competing Interests}
The authors declare that they have no relevant financial or non-financial interests to disclose.

\subsection*{Author Contributions}
MC modified the OS Terrain 50 dataset and developed the \texttt{nevis} module, with contributions from YW.
YW developed the benchmarking framework and local optima algorithm, and performed the formal analysis, with supervision provided by GRM and MC.
GRM conceptualised the study, and provided project administration and funding.  
All authors contributed to validation, reviewed and edited the manuscript, and approved the final version.

\subsection*{Data Availability}
The OS Terrain 50 data set is made openly available by Great Britain's Ordnance Survey under an Open Government Licence\footnote{\url{https://www.nationalarchives.gov.uk/doc/open-government-licence/version/3/}}.
Although these data are not redistributed as part of this paper or the accompanying software, we acknowledge its use by including the requested statement `Contains OS data \copyright\ Crown Copyright [and database right] (2024)'.

The \texttt{nevis} module created for this project that downloads and modifies the OS Terrain 50 data set can be downloaded from GitHub (\url{https://github.com/CardiacModelling/BenNevis}) or PyPI\footnote{\url{https://pypi.org/project/nevis/}}.

The code needed to re-run the experiments and generate the figures shown in this study is available on \url{https://github.com/CardiacModelling/BenNevisBenchmark}. 
Sample animations for illustrating algorithm behaviour can also be found in the same repository.

Permanent archives of code are available on Zenodo at \url{https://doi.org/10.5281/zenodo.13883484} for the \texttt{nevis} package and \url{https://doi.org/10.5281/zenodo.13883478} for the benchmarking code.

\subsection*{Open Access}
This research was funded in whole, or in part, by the Wellcome Trust [212203/Z/18/Z]. 
For the purpose of open access, the authors have applied a CC-BY public copyright licence to any Author Accepted Manuscript version arising from this submission.

\clearpage
\appendix

\section{Defining the sea bed}\label{app:sea}

First, we defined a `mask', assigning every data point a label of either `sea' or `land'.
To generate this mask, we iterated over the data in blocks of 80 by 80 points, starting with the top (most northern) row of blocks, moving left-to-right (west-to-east), then moving south, and continuing in an inward-spiralling fashion.
Within each block, the points were assigned a `sea' label if they had a height below 0\,m (or no height was available) and were adjacent to either the map edge or another point with a `sea' designation.
This was repeated until the labels within the block no longer changed, before moving on to the next block.
The entire process was then repeated (again spiralling over the blocks from top-left to top-right and onwards) until the labels had all stabilised.

As a result of this strategy, bays and the mouths of rivers were classified as `sea', although in all but two cases this stopped a few km inland, as the river became narrower than the 50\,m data resolution.
Three data points were manually elevated to a height of 0.01m to stop dry but below sea-level areas of land (e.g. Holme Fen) being classified as sea: one in the Great Ouse river (Norfolk, Ordnance Survey National Grid square TF 50), one in the river Yare (TG 50), and one in Oulton Dyke (TM59).
Seven OS Terrain 50 `tiles' corresponding to coastal areas (grid squares NT68, NR56, NR57, NR35, NR24, NR34, NR44, and NR33) contained points with elevation 0\,m which were misclassified as land: these were manually lowered to below zero.
After these manual interventions, the sea-mask algorithm was re-run.

We started by labelling all points in the sea mask as `untreated'.
Next, we created a set, $S_1$, of all `untreated' points adjacent to land.
We then iteratively:
(1) set the height of all points in $S_i$ to 0.01m below their lowest neighbour not in the untreated set; 
(2) Marked all points in $S_i$ as `treated`; 
(3) Created a set $S_{i+1}$ of all untreated `sea' points neighbouring $S_i$; 
and (4) incremented $i$ and repeated the steps until $S_{i+1}$ became an empty set.
This resulted in the progressively deepening sea visible in \fig{map} (dark blue colours).
The final data set does not contain the `sea/land' or `treated/untreated' bit masks, but consists only of the grid of elevation points.

\section{Further performance measures}
\label{app:further-pm}
We now list some further performance measures for optimisers. 
We use the same variables as in Section~\ref{sec:performance}  and denote \(N_{us} = N - N_s\) for the number of unsuccessful runs.

\begin{definition}
	The success performance \citep{suganthan_problem_2005} is defined as
	\[
		\text{SP} = \frac{T_s}{p_s}.
	\]
\end{definition}
A lower SP indicates better performance. 
SP is a variant of ERT where it is assumed that the mean running time \(T_{s}\) of successful runs is equal to the mean running time of all runs.

\begin{definition}
	The penalised average runtime \citep{kerschke2018parameterization} is defined as
	\[
		\text{PAR} = \frac 1 N \left(k N_{us} T_\text{max} + \sum_{i \in I_s} T_i \right),
	\]
	where \(k\) is a constant called the \textit{penalty factor}.
\end{definition}
A lower \(\text{PAR}\) indicates better performance.
PAR is the mean running time of all runs but with the unsuccessful ones penalised by a factor \(k\).
Conventionally, the penalty factor \(k\) is taken to be \(2\) or \(10\), which gives two performance measures \(\text{PAR2}\) and \(\text{PAR10}\).

\begin{definition}
	The dominated hypervolume \citep{bossek_multi-objective_2020} is defined as
	\[
		\text{HV} = p_s(T_\text{max} - T_s).
	\]
\end{definition}
A higher HV indicates better performance.

We mention that the proportion of dark blue (Ben Nevis) region in the height-band graph (Figure~\ref{fig:hb-1}) is proportional to the performance measure HV. 
Indeed, the area \(A\) of this region can be calculated as

\begin{align*}
	A & = \sum_{j=1}^{T_\text{max}} \sum_{i=1}^n [F_{i,j} \geqslant f_\text{target}]
	  =  \sum_{i=1}^n \sum_{j=1}^{T_\text{max}} [F_{i,j} \geqslant f_\text{target}] 
	   =  \sum_{i\in I_s} \left( T_\text{max} - T_i \right)                   
	   = N_s (T_\text{max} - T_s)                                            
	   = N (\text{HV}),
\end{align*}
where \([p] = 1\) if boolean expression \(p\) is true and \(0\) otherwise.
Then the proportion of this region in the whole graph is \(A / (N T_\text{max}) = \text{HV} / T_\text{max}\).

\end{document}